# A QUASI-ADDITIVE PROPERTY OF HOMOLOGICAL SHIFT IDEALS

SHAMILA BAYATI

Abstract. In this paper, we investigate which classes of monomial ideals have a quasi-additive property of homological shift ideals. More precisely, for a monomial ideal $I$ we are interested to find out whether $\mathrm{HS}_{i+j}(I) \subseteq \mathrm{HS}_i(\mathrm{HS}_j(I))$. It turns out that $\mathbf{c}$-bounded principal Borel ideals as well as polymatroidal ideals satisfying strong exchange property, and polymatroidal ideals generated in degree two have this quasi-additive property. For squarefree Borel ideals, we even have equality. Besides, the inclusion holds for every equigenerated Borel ideal and polymatroidal ideal when $j = 1$.

## Introduction

A recent approach in studying syzygies of a multigraded module is considering the ideals generated by their multigraded shifts which following [9] we call them *homological shift ideals*. It first came up during a discussion among Jürgen Herzog, Somayeh Bandari, and the author in 2012 whether the property of being polymatroidal is inherited by homological shift ideals. Later it turned out that this question has a positive answer for matroidal ideals [1], polymatroidal ideals with strong exchange property [9], and polymatroidal ideals generated in degree two [7]. Besides, other properties inherited by homological shift ideals, like being (squarefree) Borel or having linear quotients, are studied in [2] and [9]. In this paper, we are mainly going to discuss a property of homological shift ideals which we call it quasi-additive property.

To be more precise, let $S = k[x_1, \dots, x_n]$ be the polynomial ring in the variables $x_1, \dots, x_n$ over a field $k$ with its natural multigrading. Throughout, a monomial and its multidegree will be used interchangeably, and $S(\mathbf{x}^{\mathbf{a}})$ will denote the free $S$-module with one generator of multidegree $\mathbf{x}^{\mathbf{a}}$. A monomial ideal $I \subseteq S$ has a (unique up to isomorphism) minimal multigraded resolution

$$\mathbf{F} : 0 \to F_p \to \dots \to F_1 \to F_0$$

with

$$F_i = \bigoplus_{\mathbf{a} \in \mathbb{Z}^n} S(\mathbf{x}^{\mathbf{a}})^{\beta_{i,\mathbf{a}}}.$$

2020 *Mathematics Subject Classification.* 13D02, 13A02, 13F20, 05E40.

The $i$th homological shift ideal of $I$ denoted by $\mathrm{HS}_i(I)$ is the ideal generated by the $i$th multigraded shifts of $I$, that is,

$$\mathrm{HS}_i(I) = (\{\mathbf{x}^{\mathbf{a}} | \ \beta_{i,\mathbf{a}} \neq 0\}).$$

Along with other results, Herzog et al. show in [9, Proposition 1.4] that if $I$ has linear quotients, then

$$\mathrm{HS}_{i+1}(I) \subseteq \mathrm{HS}_1(\mathrm{HS}_i(I))$$

for all $i$. Later it is shown in [10, Corollary 4.2] and in [5, Proposition 2.4] that if $I$ is an equigenerated squarefree Borel ideal or a matroidal ideal, then one has

$$\mathrm{HS}_{i+1}(I) = \mathrm{HS}_1(\mathrm{HS}_i(I))$$

for all $i$. So the following question naturally arises that for which classes of monomial ideals one has

$$\mathrm{HS}_{i+j}(I) \subseteq \mathrm{HS}_i(\mathrm{HS}_j(I))$$

for all $i, j$. We say that $I$ has the *quasi-additive property for homological shift ideals* or simply $I$ is *quasi-additive* if the above question has a positive answer for $I$.

In this paper, we are about to find classes of quasi-additive ideals. We first show in Theorem 1.2 that when $I$ is an equigenerated monomial ideal, $<$ is a monomial order which extends $x_1 > x_2 > \ldots > x_n$, and $I$ and $\mathrm{HS}_j(I)$ have linear quotients with respect to $<$ for some $j$, then $\mathrm{HS}_{i+j}(I) \subseteq \mathrm{HS}_i(\mathrm{HS}_j(I))$ for all $i$. This implies that $\mathbf{c}$-bounded principal Borel ideals, polymatroidal ideals satisfying strong exchange property, and the edge ideal of the complement of path graphs are among the quasi-additive ideals. It is shown in [10, Corollary 4.2], if $I$ is an equigenerated squarefree Borel ideal then $\mathrm{HS}_{i+j}(I) = \mathrm{HS}_i(\mathrm{HS}_j(I))$ for all $i, j$. We generalize this result for (not necessarily equigenerated) squarefree Borel ideals in Theorem 1.9.

In Theorem 2.1, we will show that the adjacency ideal of a polymatroidal ideal is polymatroidal as well. This, in particular, implies that the first homological shift ideal of a polymatroidal ideal is also polymatroidal; a result that has been proved by Ficarra by a different approach in [5]. So, as stated in Corollary 2.3, when $I$ is a polymatroidal ideal, one has $\mathrm{HS}_{i+1}(I) \subseteq \mathrm{HS}_i(\mathrm{HS}_1(I))$ for each $i$.

We call a monomial $\mathbf{x}^{\mathbf{a}} \in k[x_1, \ldots, x_n]$ quasi-squarefree if $\mathbf{a}$ is componentwise less than or equal to $\mathbf{1} + \hat{i}$ for some $i$ where $\mathbf{1} = (1, \ldots, 1) \in \mathbb{Z}^n$, and $\hat{i}$ is the $i$th canonical basis vector of $\mathbb{R}^n$. If $I \subseteq S$ is a monomial ideal, we define an operation that assigns to $I$ its quasi-squarefree part which is the monomial ideal generated by quasi-squarefree monomials in $\mathrm{G}(I)$. We first show in Lemma 2.5 that if we start with a polymatroidal ideal generated by quasi-squarefree monomials, then quasi-squarefree part of its adjacency ideal is also polymatroidal. Next, it turns out in Lemma 2.7 that when $I$ is a polymatroidal ideal generated in degree two, the ideal $\mathrm{HS}_i(I)$ can be obtained by taking $i$ times iterated adjacency ideals and then quasi-squarefree part, one after another, starting from $I$. On the one hand, this implies the quasi-additive property for homological shift ideals of polymatroidal ideals generated in degree two, as one can see in Theorem 2.8. On the other hand, as a result, one obtains a very recent result by Ficarra and Herzog which gives a positive answer to the conjecture about homological shift ideals of polymatroidal ideals when we restrict ourselves to those generated in degree two; see Corollary 2.9. Finally, in Proposition 2.10 via the concept of adjacency ideals, we prove $\mathrm{HS}_{i+j}(I) = \mathrm{HS}_i(\mathrm{HS}_j(I))$ when $I$ is a matroidal ideal, as one has by [5, Proposition 2.4].

## 1. Quasi-additive property for Borel ideals

Throughout, $S = k[x_1, \dots, x_n]$ denotes a polynomial ring over a field $k$ with its natural multigrading. Moreover, a monomial $\mathbf{x}^{\mathbf{a}} = x_1^{a_1} \dots x_n^{a_n}$ and its multidegree $(a_1, \dots a_n)$ will be used interchangeably. Besides, in the case that $\mathbf{x}^{\mathbf{a}}$ is a squarefree monomial, we may use its support instead of it. So we will apply some notions related to monomials (resp. squarefree monomials) for vectors in $\mathbb{Z}^n_{\geq 0}$ (resp. the subsets of [n]). If $u, v \in S$ are monomials, then $u : v$ denotes the monomial $\frac{u}{\gcd(u,v)}$. For a monomial $u \in S$, we set $\max u = \max\{k \colon x_k \text{ divides } u\}$. When $\ell = \max u$, we may sometimes write $x_\ell = \max u$ for ease of use.

Let $I \subseteq S$ be a monomial ideal. We denote its minimal set of monomial generators by $G(I)$. A monomial ideal $I \subseteq S$ is said to have linear quotients if there exists an ordering $u_1, \dots, u_r$ of the elements of $\mathrm{G}(I)$ such that for each $i = 1, \dots, r-1$, the colon ideal $(u_1, \dots, u_i) : (u_{i+1})$ is generated by a subset of $\{x_1, \dots, x_n\}$. If $I$ has linear quotients with respect to the ordering $u_1, \dots, u_r$ of $\mathrm{G}(I)$, then

$$\{x_j : x_j \in (u_1, \dots, u_i) : (u_{i+1})\}$$

is denoted by $\mathrm{set}(u_{i+1})$.

**Remark 1.1.** Let a monomial ideal $I \subseteq S$ have linear quotients. By [11, Lemma 1.5], a minimal multigraded free resolution $\mathbf{F}$ of $I$ can be described as follows: the $S$-module $F_i$ in homological degree $i$ of $\mathbf{F}$ is the multigraded free $S$-module whose basis is formed by monomials $ux_{\ell_1} \dots x_{\ell_i}$ which $u \in \mathrm{G}(I)$ and $x_{\ell_1}, \dots, x_{\ell_i}$ are distinct elements of $\mathrm{set}(u)$.

**Theorem 1.2.** *Let $I$ be a monomial ideal generated in a single degree and let $j$ be a nonnegative integer. Suppose that $<$ is a monomial order which extends $x_1 > x_2 > \dots > x_n$. If the ideals $I$ and $\mathrm{HS}_j(I)$ have linear quotients with respect to the descending order of their minimal set of monomial generators by $<$, then for every $i$*

$$\mathrm{HS}_{i+j}(I) \subseteq \mathrm{HS}_i(\mathrm{HS}_j(I)).$$

*Proof.* We show that each generator $ux_{\ell_1} \dots x_{\ell_{i+j}}$ of $\mathrm{HS}_{i+j}(I)$ with $u \in \mathrm{G}(I)$ and $\{\ell_1 < \dots < \ell_{i+j}\} \subseteq \mathrm{set}(u)$ belongs to $\mathrm{HS}_i(\mathrm{HS}_j(I))$. Notice that by Remark 1.1, $w = ux_{\ell_{i+1}} \dots x_{\ell_{i+j}} \in \mathrm{HS}_j(I)$. Besides, for each $t = 1, \dots, i$ one has

$$w_t = ux_{\ell_t}\widehat{x_{\ell_{i+1}}}x_{\ell_{i+2}} \dots x_{\ell_{i+j}} = ux_{\ell_t}x_{\ell_{i+2}} \dots x_{\ell_{i+j}} \in \mathrm{HS}_j(I),$$

where $\widehat{x_{\ell_{i+1}}}$ denotes omitted variable in the product. Moreover, $w$ and $w_t$'s belong to the minimal set of monomial generators of $\mathrm{HS}_j(I)$ because $I$ is generated in a single degree. Since $x_{\ell_t} > x_{\ell_{i+1}}$ for each $t = 1, \dots, i$ by assumption, multiplying this inequality by $ux_{\ell_{i+2}} \dots x_{\ell_{i+j}}$ yields that $w_t > w$. In addition,

$$w_t \colon w = x_{\ell_t}.$$

Hence with respect to the descending order of the minimal set of monomial generators of $\mathrm{HS}_j(I)$ by $<$, one has

$$x_{\ell_t} \in \mathrm{set}(w)$$

for each $t = 1, \dots, i$. In particular, by Remark 1.1,

$$ux_{\ell_1} \dots x_{\ell_{i+j}} = wx_{\ell_1} \dots x_{\ell_i} \in \mathrm{HS}_i(\mathrm{HS}_j(I)),$$

as desired. □

Let $\mathbf{c}$ be a vector in $\mathbb{Z}^n$ with non-negative entries. A monomial $\mathbf{x}^{\mathbf{b}} \in S$ is called $\mathbf{c}$-bounded if $\mathbf{b}$ is componentwise less than or equal to $\mathbf{c}$. Associated to each monomial ideal $I \subseteq S$, $I^{\leq \mathbf{c}}$ denotes the monomial ideal

$$I^{\leq \mathbf{c}} = (\mathbf{x}^{\mathbf{b}} \colon \mathbf{x}^{\mathbf{b}} \in \mathrm{G}(I) \text{ and } \mathbf{x}^{\mathbf{b}} \text{ is } \mathbf{c}\text{-bounded}) \subseteq S.$$

The ideal $I$ is called $\mathbf{c}$-bounded if $I = I^{\leq \mathbf{c}}$. Notice that each squarefree monomial ideal is $\mathbf{c}$-bounded for $\mathbf{c} = (1, 1, \dots, 1)$.

An operation that sends a monomial $u$ to $(u/x_j)x_i$ is called a *Borel move* if $x_j$ divides $u$ and $i < j$. When $u$ is a $\mathbf{c}$-bounded (resp. squarefree) monomial, a Borel move is called a $\mathbf{c}$*-bounded (resp. squarefree) Borel move* if the monomial $(u/x_j)x_i$ is also $\mathbf{c}$-bounded (resp. squarefree). A monomial ideal $I \subseteq S$ is called a Borel ideal if it is closed under Borel moves. The ideal $I$ is called $\mathbf{c}$-bounded (resp. squarefree) Borel, if it is a $\mathbf{c}$-bounded (resp. squarefree) monomial ideal and closed under $\mathbf{c}$-bounded (resp. squarefree) Borel moves. A subset $B$ of a Borel ideal $I$ is called its *Borel generator* if $I$ is the smallest Borel ideal containing $B$. A Borel ideal $I$ is called a *principal Borel ideal* if it has a Borel generator of cardinality one.

**Corollary 1.3.** *Let $I$ be a $\mathbf{c}$-bounded principal Borel ideal. Then*

$$\mathrm{HS}_{i+j}(I) \subseteq \mathrm{HS}_j(\mathrm{HS}_i(I))$$

*for each $i, j$.*

*Proof.* By [9, Theorem 2.2], if $I$ is a $\mathbf{c}$-bounded principal Borel ideal, then $\mathrm{HS}_j(I)$ has linear quotients for each $j$. Indeed by proof of [2, Theorem 2.4] and [9, Proposition 2.6], it turns out that each ideal $\mathrm{HS}_j(I)$ has linear quotients when the elements of $\mathrm{HS}_j(I)$ are ordered decreasingly with respect to the lexicographical order with $x_1 > x_2 > \dots > x_n$, as required in Theorem 1.2. Hence for every $i$

$$\mathrm{HS}_{i+j}(I) \subseteq \mathrm{HS}_i(\mathrm{HS}_j(I)).$$

□

**Example 1.4.** Consider the principal Borel ideal $I \subseteq k[x_1, x_2, x_3]$ with Borel generator $\{x_1x_2x_3\}$, that is,

$$I = (x_1^3, x_1^2x_2, x_1^2x_3, x_1x_2^2, x_1x_2x_3).$$

Then one has

$$\begin{aligned} \mathrm{HS}_1(I) &= (x_1^3x_2, x_1^3x_3, x_1^2x_2^2, x_1^2x_2x_3, x_1x_2^2x_3); \\ \mathrm{HS}_2(I) &= (x_1^3x_2x_3, x_1^2x_2^2x_3). \end{aligned}$$

Besides, $\mathrm{HS}_1(\mathrm{HS}_1(I)) = (x_1^3x_2^2, x_1^3x_2x_3, x_1^2x_2^2x_3)$. Recall that the ideal $I$, and by [2, Theorem 2.4] the ideal $\mathrm{HS}_1(I)$ have linear quotients with respect to the lexicographical order induced by $x_1 > x_2 > \ldots > x_n$. Hence this example shows that equality does not necessarily hold in Theorem 1.2.

**Corollary 1.5.** *Let $I$ be an equigenerated Borel ideal. Then*

$$\mathrm{HS}_{i+1}(I) \subseteq \mathrm{HS}_i(\mathrm{HS}_1(I))$$

*for each $i$.*

*Proof.* By [2, Proposition 2.2], the ideal $\mathrm{HS}_1(I)$ has linear quotients with respect to the lexicographical order induced by the ordering $x_1 > x_2 > \ldots > x_n$ of variables. Now the assertion follows from Theorem 1.2. □

**Remark 1.6.** Let $I$ be the edge ideal of the complement of a path graph. By [9, Proposition 4.2], for each $j$ the ideal $\mathrm{HS}_j(I)$ has linear quotients with respect to the lexicographical order induced by $x_1 > x_2 > \ldots > x_n$. Hence by Theorem 1.2 such an ideal $I$ is quasi-additive.

**Remark 1.7.** Let $I$ be a squarefree Borel ideal. It is shown in [2, Theorem 3.3] that the ideal $\mathrm{HS}_i(I)$ has linear quotients for each $i$ with respect to the following order $w_1, \ldots, w_r$ of the minimal set of monomial generators of $\mathrm{HS}_i(I)$: $i < j$ implies that either (i) $\deg(w_i) < \deg(w_j)$ or (ii) $\deg(w_i) = \deg(w_j)$ and $w_i >_{lex} w_j$. Here lexicographical order is induced by the ordering $x_1 > x_2 > \ldots > x_n$.

**Remark 1.8.** Let $I$ be a squarefree Borel ideal. Applying [6, Theorem 2.1] and [8, Lemma 4.4.1] to the minimal multigraded free resolution described for Borel ideals in [4, Theorem 2.1], one obtains the minimal multigraded free resolution $\mathbf{F}$ of $I$ as follows: the basis of the multigraded free $S$-module $F_i$ in homological degree $i$ of $\mathbf{F}$ is formed by those multihomogeneous elements of multidegree $\mathbf{a}$ such that $\mathbf{x}^{\mathbf{a}}$ is a *squarefree* monomial $ux_{\ell_1}\ldots x_{\ell_i}$ with $u \in \mathrm{G}(I)$ and $\ell_t < \max u$ for each $t = 1, \ldots, i$. A sequence $x_{\ell_1}, \ldots, x_{\ell_i}$ satisfying these conditions is called an admissible sequence for $u$.

By [10, Corollary 4.2] if $I$ is an equigenerated squarefree Borel ideal, then one has $\mathrm{HS}_{i+j}(I) = \mathrm{HS}_i(\mathrm{HS}_j(I))$ for all $i, j$. The following result gives a generalization for (not necessarily equigenerated) squarefree Borel ideals.

**Theorem 1.9.** *Let $I$ be a squarefree Borel ideal. Then*

$$\mathrm{HS}_{i+j}(I) = \mathrm{HS}_i(\mathrm{HS}_j(I)).$$

*for each $i, j$.*

*Proof.* The assertion is trivial if $j = 0$. So assume that $j > 0$. We first show that $\mathrm{HS}_{i+j}(I) \subseteq \mathrm{HS}_i(\mathrm{HS}_j(I))$. Recall the description of the minimal multigraded free resolution of $I$ in Remark 1.8. Let $ux_{\ell_1}\ldots x_{\ell_{i+j}} \in \mathrm{HS}_{i+j}(I)$ where $u \in \mathrm{G}(I)$ and $x_{\ell_1}, \ldots, x_{\ell_{i+j}}$ is an admissible sequence for $u$ with $\ell_1 < \cdots < \ell_{i+j}$. One also has

$$ux_{\ell_{i+1}}\ldots x_{\ell_{i+j}} \in \mathrm{HS}_j(I),$$

however, this monomial may not belong to the minimal set of monomial generators of $\mathrm{HS}_j(I)$. Assume that $w = vx_{k_1}\dots x_{k_j}$ is a squarefree monomial in $\mathrm{G}(\mathrm{HS}_j(I))$ which divides $ux_{\ell_{i+1}}\dots x_{\ell_{i+j}}$. Here $v \in \mathrm{G}(I)$ and $x_{k_1},\dots,x_{k_j}$ is an admissible sequence for $v$ with $k_1 < \cdots < k_j$. Now recall that $\mathrm{HS}_j(I)$ has linear quotients as clarified in Remark 1.7. So it is enough to show that $x_{\ell_t} \in \mathrm{set}(w)$ for each $t = 1,\dots,i$ which by Remark 1.1 implies that

$$wx_{\ell_1}\dots x_{\ell_i} = (vx_{k_1}\dots x_{k_j})x_{\ell_1}\dots x_{\ell_i} \in \mathrm{HS}_i(\mathrm{HS}_j(I)).$$

Consequently, since this monomial divides $ux_{\ell_1}\dots x_{\ell_{i+j}}$, we will obtain that

$$ux_{\ell_1}\dots x_{\ell_{i+j}} \in \mathrm{HS}_i(\mathrm{HS}_j(I),$$

as desired. One has $\ell_t \neq k_j$ for each $t = 1,\dots,i$. So two cases may happen for each $t = 1,\dots,i$:

Case 1. If $\ell_t < k_j$, we set

$$w_t = vx_{\ell_t}x_{k_1}\dots\widehat{x_{k_j}},$$

where $\widehat{x_{k_j}}$ denotes an omitted variable in the product. It is clear that $x_{\ell_t}, x_{k_1},\dots,x_{k_{j-1}},\widehat{x_{k_j}}$ is an admissible sequence for $v$. So $w_t = vx_{\ell_t}x_{k_1}\dots\widehat{x_{k_j}} \in \mathrm{HS}_j(I)$. Suppose that $\tilde{w}_t$ is an element of $\mathrm{G}(\mathrm{HS}_j(I))$ that divides $w_t$. Since $w$ is also an element of $\mathrm{G}(\mathrm{HS}_j(I))$, we conclude that $x_{\ell_t}$ must divide $\tilde{w}_t$. On the other hand, $\tilde{w}_t$ comes before $w$ in the order of generators of $\mathrm{HS}_j(I)$ described in Remark 1.7. Thus $\tilde{w}_t : w = x_{\ell_t} \in \mathrm{set}(w)$.

Case 2. If $k_j < \ell_t$, we set

$$w_t = (\frac{v}{\max v}x_{\ell_t})x_{k_1}\dots x_{k_j}.$$

The condition $k_j < \ell_t$ implies that $u \neq v$. So we may assume that

(1) $$\deg(v) < \deg(u).$$

To prove $\frac{v}{\max v}x_{\ell_t} \in I$, we claim that at least one of the variables $x_{\ell_{i+1}},\dots,x_{\ell_{i+j}}$ divides $v$, say $x_{\ell_s}$ which implies that $\ell_t < i+1 \leq \ell_s \leq \max v$ and consequently $\frac{v}{\max v}x_{\ell_t} \in I$. Assume on the contrary that none of the variables $x_{\ell_{i+1}},\dots,x_{\ell_{i+j}}$ divide $v$. Since, on the other hand, $v$ divides $w$, and $w$ divides the squarefree monomial $ux_{\ell_{i+1}}\dots x_{\ell_{i+j}}$, we deduce that $v$ divides $u$; a contradiction to the fact that by (1) $u$ and $v$ are distinct elements of $\mathrm{G}(I)$. Next, notice that the assumption $k_j < \ell_t$ guarantees that $x_{k_1},\dots,x_{k_j}$ with $k_1 < \cdots < k_j$ is an admissible sequence for $\frac{v}{\max v}x_{\ell_t}$.
Considering an element $\tilde{w}_t \in \mathrm{G}(\mathrm{HS}_j(I))$ that divides $w_t$, the same argument as used in Case 1 shows that $\tilde{w}_t : w = x_{\ell_t} \in \mathrm{set}(w)$.

To finish the proof, we show the other inclusion, that is,

$$\mathrm{HS}_i(\mathrm{HS}_j(I)) \subseteq \mathrm{HS}_{i+j}(I).$$

Regarding Remark 1.7 $\mathrm{HS}_j(I)$ has linear quotients. So recall the description of generators of $\mathrm{HS}_i(\mathrm{HS}_j(I))$ by Remark 1.1 and the description of generators of $\mathrm{HS}_j(I)$ by Remark 1.8. Now suppose that the squarefree monomial

$$ux_{\ell_1}\dots x_{\ell_j}x_{k_1}\dots x_{k_i}$$

belongs to $\mathrm{HS}_i(\mathrm{HS}_j(I))$ with $x_{k_1},\ldots,x_{k_i} \in \mathrm{set}(ux_{\ell_1}\ldots x_{\ell_j})$ in the ideal $\mathrm{HS}_j(I)$, and $x_{\ell_1},\ldots,x_{\ell_j}$ is an admissible sequence for $u \in \mathrm{G}(I)$. In particular, assume that $ux_{\ell_1}\ldots x_{\ell_j}$ belongs to $\mathrm{G}(\mathrm{HS}_j(I))$. We need to show that $k_t < \max u$ for each $t = 1,\ldots,i$ to deduce that

$$x_{\ell_1},\ldots,x_{\ell_j},x_{k_1},\ldots,x_{k_i}$$

is an admissible sequence for $u$ and consequently,

$$ux_{\ell_1}\ldots x_{\ell_j}x_{k_1}\ldots x_{k_i} \in \mathrm{HS}_{i+j}(I).$$

Fix $t = 1,\ldots,i$. We have $x_{k_t} \in \mathrm{set}(ux_{\ell_1}\ldots x_{\ell_j})$ in the ideal $\mathrm{HS}_j(I)$. So there exists a squarefree monomial $vx_{s_1}\ldots x_{s_j} \in \mathrm{G}(\mathrm{HS}_j(I))$ with $v \in \mathrm{G}(I)$ and admissible sequence $x_{s_1},\ldots,x_{s_j}$ for $v$ with $s_1 < \cdots < s_j$ such that

$$vx_{s_1}\ldots x_{s_j} : ux_{\ell_1}\ldots x_{\ell_j} = x_{k_t}, \tag{2}$$

and $vx_{s_1}\ldots x_{s_j}$ comes before $ux_{\ell_1}\ldots x_{\ell_j}$ in the ordering of generators of $\mathrm{HS}_j(I)$ described in Remark 1.7. Thus one has either

$$\deg(vx_{s_1}\ldots x_{s_j}) < \deg(ux_{\ell_1}\ldots x_{\ell_j}) \tag{3}$$

or $\deg(vx_{s_1}\ldots x_{s_j}) = \deg(ux_{\ell_1}\ldots x_{\ell_j})$ and $vx_{s_1}\ldots x_{s_j} >_{lex} ux_{\ell_1}\ldots x_{\ell_j}$.

First, assume that $\deg(vx_{s_1}\ldots x_{s_j}) < \deg(ux_{\ell_1}\ldots x_{\ell_j})$. Regarding (2), since we are working with squarefree monomials, we conclude that $k_t \neq \max u$. On the contrary, suppose that $k_t > \max u$. Thus

$$\max(vx_{s_1}\ldots x_{s_j}) = \max v \geq k_t > \max u = \max(ux_{\ell_1}\ldots x_{\ell_j}).$$

As a result $\max(vx_{s_1}\ldots x_{s_j}) = \max v$ does not divide $ux_{\ell_1}\ldots x_{\ell_j}$. So regarding (2),

$$\max v = k_t. \tag{4}$$

Set

$$p = \max\{r\colon\ \ x_r | ux_{\ell_1}\ldots x_{\ell_j}\ \text{ and }\ x_r \not| vx_{s_1}\ldots x_{s_j}\}.$$

Consider the admissible sequence $x_{s_1},\ldots,x_{s_{j-1}},x_p$ for $(v/\max v)x_{s_j}$ when $p < s_j$. Furthermore, regarding $p \leq \max u < k_t = \max v$ consider the element $(v/\max v)x_p$ with the admissible sequence $x_{s_1},\ldots,x_{s_{j-1}},x_{s_j}$ when $s_j < p$. Both admissible sequences give an element

$$(v/\max v)x_{s_1}\ldots x_{s_{j-1}}x_{s_j}x_p$$

of the ideal $\mathrm{HS}_j(I)$. By (2) and (4), the monomial $(v/\max v)x_{s_1}\ldots x_{s_{j-1}}x_{s_j}x_p \in \mathrm{HS}_j(I)$ with the same degree as $vx_{s_1}\ldots x_{s_j}$ divides

$$ux_{\ell_1}\ldots x_{\ell_j} \in \mathrm{G}(\mathrm{HS}_j(I)),$$

a contradiction to (3). Hence in the case that $\deg(vx_{s_1}\ldots x_{s_j}) < \deg(ux_{\ell_1}\ldots x_{\ell_j})$, we have $k_t < \max u$, as desired.

Next assume that

$$\deg(vx_{s_1}\ldots x_{s_j}) = \deg(ux_{\ell_1}\ldots x_{\ell_j})\ \text{ and }\ vx_{s_1}\ldots x_{s_j} >_{lex} ux_{\ell_1}\ldots x_{\ell_j}.$$

Then (2) along with the lexicographical order of generators immediately yields that $k_t < \max(ux_{\ell_1}\ldots x_{\ell_j}) = \max u$. □

## 2. Quasi-additive property for polymatroidal ideals

In this section, we consider polymatroidal ideals and study the quasi-additive property for some important classes of these ideals.

Let $\mathbf{a}, \mathbf{b} \in \mathbb{Z}^n_{\geq 0}$ be two vectors. The *join* $\mathbf{a} \vee \mathbf{b}$ of $\mathbf{a}$ and $\mathbf{b}$ denotes the vector in $\mathbb{Z}^n_{\geq 0}$ corresponding to $\operatorname{lcm}(\mathbf{x^a}, \mathbf{x^b})$. For each $i \in [n]$ we will denote the $i$th canonical basis vector of $\mathbb{R}^n$ by $\hat{i}$. We consider the distance between the monomials of $S$ in the sense of [3], that is,

$$d(\mathbf{x^a}, \mathbf{x^b}) = \frac{1}{2} \sum_{k=1}^{n} |\deg_k \mathbf{x^a} - \deg_k \mathbf{x^b}|$$

where for a monomial $\mathbf{x^a} = x_1^{a_1} \dots x_n^{a_n}$, one has $\deg_k \mathbf{x^a} = a_k$. Following [1], we consider the *adjacency graph* $G_I$ of $I$ whose set of vertices is the set of minimal monomial generators $\mathrm{G}(I)$ of $I$, and two vertices $\mathbf{x^a}$ and $\mathbf{x^b}$ are adjacent if $d(\mathbf{x^a}, \mathbf{x^b}) = 1$. The *adjacency ideal* of $I$, denoted by $\mathrm{A}(I)$, is defined to be the monomial ideal generated by the least common multiples of adjacent vertices in $G_I$, that is,

$$\mathrm{A}(I) = \langle \operatorname{lcm}(\mathbf{x^a}, \mathbf{x^b}) : d(\mathbf{x^a}, \mathbf{x^b}) = 1 \rangle \subseteq k[x_1, \dots, x_n].$$

In terms of bases, a discrete polymatroid $\mathcal{P}$ is a pair $([n], \mathcal{B})$ which the nonempty finite set of bases $\mathcal{B} \subseteq \mathbb{Z}^n_{\geq 0}$ is satisfying the following conditions:

(I) Every $\mathbf{a}, \mathbf{b} \in \mathcal{B}$ have the same modulus, that is, $a_1 + \dots + a_n = b_1 + \dots + b_n$;

(II) If $\mathbf{a}, \mathbf{b} \in \mathcal{B}$, for each $i$ with $a_i > b_i$, there exists $j \in [n]$ such that $b_j > a_j$ and $\mathbf{a} - \hat{i} + \hat{j} \in \mathcal{B}$.

Property (II) is called the *exchange property*. It is known that bases of $\mathcal{P}$ possess the following *symmetric exchange property*:

> If $\mathbf{a}, \mathbf{b} \in \mathcal{B}$, for each $i$ with $a_i > b_i$, there exists $j \in [n]$ such that $b_j > a_j$ and $\mathbf{a} - \hat{i} + \hat{j}, \mathbf{b} - \hat{j} + \hat{i} \in \mathcal{B}$.

The discrete polymatroid $\mathcal{P} = ([n], \mathcal{B})$ is said to satisfy the *strong exchange property* if for every $\mathbf{a}, \mathbf{b} \in \mathcal{B}$ and each $i, j \in [n]$ with $a_i > b_i$ and $a_j < b_j$, one has $\mathbf{a} - \hat{i} + \hat{j} \in \mathcal{B}$. A monomial ideal $I \subseteq S$ is called a *polymatroidal ideal* if its minimal set of monomial generators $\mathrm{G}(I)$ corresponds to the set of bases of a discrete polymatroid. A *matroidal ideal* is a squarefree polymatroidal ideal.

**Theorem 2.1.** *Let $I$ be a polymatroidal ideal. Then its adjacency ideal* $\mathrm{A}(I)$ *is also a polymatroidal ideal.*

*Proof.* Let $\mathcal{P} = ([n], \mathcal{B})$ be the discrete polymatroid corresponding to $I$. For each basis $\mathbf{a}$ of $\mathcal{P}$ we define $\operatorname{set}(\mathbf{a}) \subseteq [n]$ as follows: $i \in \operatorname{set}(\mathbf{a})$ if there exists a vector $\mathbf{b} \in \mathcal{B}$ and $j \neq i$ in $[n]$ such that $\mathbf{a} + \hat{i} = \mathbf{b} + \hat{j}$.

Remind that we use monomials in $k[x_1, \dots, x_n]$ and vectors in $\mathbb{Z}^n_{\geq 0}$ interchangeably. Consider two distinct arbitrary elements $\mathbf{b} + \hat{i}, \mathbf{c} + \hat{j} \in \mathrm{A}(I)$ with $\mathbf{b}, \mathbf{c} \in \mathcal{B}$ ,

$i \in \text{set}(\mathbf{b})$, and $j \in \text{set}(\mathbf{c})$. Since $i \in \text{set}(\mathbf{b})$, there exists $\mathbf{b}' \in \mathcal{B}$ and $\ell \in [n]$ such that $\deg_\ell \mathbf{b} > \deg_\ell \mathbf{b}'$, and

$$\mathbf{b} + \hat{\imath} = \mathbf{b}' + \hat{\ell} \tag{5}$$

We are going to check the exchange property for the elements $\mathbf{b} + \hat{\imath}$ and $\mathbf{c} + \hat{\jmath}$, that is, for each element $b \in [n]$ with $\deg_b(\mathbf{b}+\hat{\imath}) > \deg_b(\mathbf{c}+\hat{\jmath})$, we find an element $c \in [n]$ such that

$$\deg_c(\mathbf{c} + \hat{\jmath}) > \deg_c(\mathbf{b} + \hat{\imath}) \text{ and } \mathbf{b} + \hat{\imath} - \hat{b} + \hat{c} \in \mathrm{A}(I). \tag{6}$$

For such an element $b$, we may assume that $\deg_b \mathbf{b} > \deg_b \mathbf{c}$. Otherwise $b = i$, and we can proceed with the other presentation, namely the presentation $\mathbf{b}' + \hat{\ell}$ given in (5). We thus assume that

$$\deg_b \mathbf{b} > \deg_b \mathbf{c}, \tag{7}$$

and then exchange property for bases of $\mathcal{B}$ yields that the following set is not empty:

$$T = \{c \in [n] : \deg_c \mathbf{c} > \deg_c \mathbf{b} \text{ and } \mathbf{b} - \hat{b} + \hat{c} \in \mathcal{B}\}.$$

There exist two cases:

*Case 1.* First, suppose that $i \in T$. So $\mathbf{b}_1 = \mathbf{b} - \hat{b} + \hat{\imath} \in \mathcal{B}$. If $\mathbf{b}_1 = \mathbf{c}$, then $j \neq b$ because $\mathbf{b} + \hat{\imath}$ and $\mathbf{c} + \hat{\jmath}$ are distinct elements. This implies that

$$\deg_j \mathbf{c} = \deg_j \mathbf{b}_1 \geq \deg_j \mathbf{b}.$$

In particular, the inequality is strict if $i = j$. Hence $\deg_j(\mathbf{c} + \hat{\jmath}) > \deg_j(\mathbf{b} + \hat{\imath})$. So the choice of $c = j$, for which $\mathbf{b}+\hat{\imath}-\hat{b}+\hat{\jmath} = \mathbf{c}+\hat{\jmath} \in \mathrm{A}(I)$, has the required property mentioned in (6). So we are done when $\mathbf{b}_1 = \mathbf{c}$. Otherwise, if $\mathbf{b}_1 \neq \mathbf{c}$, there exists $b' \in [n]$ such that $\deg_{b'} \mathbf{b}_1 > \deg_{b'} \mathbf{c}$, and consequently by exchange property, there exists an element $c \in [n]$ such that $\deg_c \mathbf{c} > \deg_c \mathbf{b}_1$ and

$$\mathbf{b}_2 = \mathbf{b}_1 - \hat{b'} + \hat{c} \in \mathcal{B}.$$

The vertices $\mathbf{b}_1$ and $\mathbf{b}_2$ are adjacent, Therefore, $\mathbf{b}_1 \vee \mathbf{b}_2 \in \mathrm{A}(I)$. But

$$\mathbf{b}_1 \vee \mathbf{b}_2 = \mathbf{b}_1 \vee (\mathbf{b}_1 - \hat{b'} + \hat{c}) = \mathbf{b}_1 + \hat{c} = (\mathbf{b} - \hat{b} + \hat{\imath}) + \hat{c} \in \mathrm{A}(I). \tag{8}$$

On the other hand, $c \neq b$ because $\deg_c \mathbf{c} > \deg_c \mathbf{b}_1$ but

$$\deg_b \mathbf{b}_1 = \deg_b \mathbf{b} - 1 \geq \deg_b \mathbf{c}$$

by (7). This implies that

$$\deg_c \mathbf{c} > \deg_c \mathbf{b}_1 \geq \deg_c \mathbf{b},$$

and the last inequality is strict if $c = i$. Therefore, $\deg_c(\mathbf{c} + \hat{\jmath}) > \deg_c(\mathbf{b} + \hat{\imath})$. Thus, regarding (8), the element $c$ is an appropriate choice for (6).

*Case 2.* Next, suppose that there exists an element $c \in T$ such that $c \neq i$. Thus $\mathbf{b}_1 = \mathbf{b} - \hat{b} + \hat{c} \in \mathcal{B}$. Moreover, we have $\deg_c(\mathbf{c} + \hat{\jmath}) > \deg_c(\mathbf{b} + \hat{\imath})$ because $c \in T$ and $c \neq i$. So it remains to show that $\mathbf{b} - \hat{b} + \hat{c} + \hat{\imath} \in \mathrm{A}(I)$ as required in (6).

If $b = \ell$ where $\ell$ is introduced in (5), then the vertex $\mathbf{b}_1 = \mathbf{b} - \hat{\ell} + \hat{c} \in \mathcal{B}$ becomes adjacent to $\mathbf{b}' = \mathbf{b} - \hat{\ell} + \hat{\imath} \in \mathcal{B}$. Hence $\mathbf{b}' \vee \mathbf{b}_1 \in \mathrm{A}(I)$. Thus

$$\mathbf{b}' \vee \mathbf{b}_1 = (\mathbf{b}_1 - \hat{c} + \hat{\imath}) \vee \mathbf{b}_1 = \mathbf{b}_1 + \hat{\imath} = (\mathbf{b} - \hat{b} + \hat{c}) + \hat{\imath} \in \mathrm{A}(I).$$

So $c$ is the desired element.

If $b \neq \ell$, then $\mathbf{b}'$ as introduced in (5) is equal to $\mathbf{b}_1 - (\hat{c} + \hat{\ell}) + (\hat{b} + \hat{\imath})$. Notice that $d(\mathbf{b}', \mathbf{b}_1) = 2$ because $\{c, \ell\}$ and $\{b, i\}$ are disjoint sets. Now by exchange property for bases of the discrete polymatroid, there exist a common neighbor $\mathbf{b}_2 \in \mathcal{B}$ of $\mathbf{b}'$ and $\mathbf{b}_1$ by $i$ pivoted in, namely $\mathbf{b}_2 = \mathbf{b}_1 - \hat{\ell} + \hat{\imath}$ or $\mathbf{b}_2 = \mathbf{b}_1 - \hat{c} + \hat{\imath}$. In any case, $\mathbf{b}_1 \vee \mathbf{b}_2 = \mathbf{b}_1 + \hat{\imath} = \mathbf{b} - \hat{b} + \hat{c} + \hat{\imath} \in \mathrm{A}(I)$, as desired. □

**Corollary 2.2.** *(*[5, Theorem 2.2]*) Let $I$ be a polymatroidal ideal. Then* $\mathrm{HS}_1(I)$ *is also a polymatroidal ideal.*

*Proof.* Regarding Remark 1.1, one has $\mathrm{HS}_1(I) = \mathrm{A}(I)$ when $I$ is a monomial ideal with linear quotients generated in a single degree. Now recall that by [11, Lemma 1.3], polymatroidal ideals have linear quotients. □

The next result partially generalizes [5, Proposition 2.4] by Ficarra about matroidal ideals. Recall that by the result of Herzog et al. in [9, Proposition 1.4] $\mathrm{HS}_{i+1}(I) \subseteq \mathrm{HS}_1(\mathrm{HS}_i(I))$ for every $i$ if $I$ has linear quotients.

**Corollary 2.3.** *Let $I$ be a polymatroidal ideal. Then*

$$\mathrm{HS}_{i+1}(I) \subseteq \mathrm{HS}_i(\mathrm{HS}_1(I))$$

*for each $i$.*

*Proof.* By [11, Lemma 1.3] polymatroidal ideals have linear quotients with respect to the reverse lexicographical order of the generators, as required in Theorem 1.2. Now the assertion immediately follows from Theorem 1.2 and Corollary 2.2. □

**Remark 2.4.** Let $I$ be a polymatroidal ideal corresponding to a discrete polymatroid satisfying strong exchange property. Then by [9, Corollary 2.6], the ideal $\mathrm{HS}_i(I)$ is polymatroidal for every $i$. This yields that

$$\mathrm{HS}_{i+j}(I) \subseteq \mathrm{HS}_i(\mathrm{HS}_j(I)).$$

for each $i, j$ by Theorem 1.2. One can see that equality does not hold necessarily. Indeed, the ideal presented in Example 1.4 is a polymatroidal ideal with strong exchange property for which $\mathrm{HS}_2(I) \subsetneq \mathrm{HS}_1(\mathrm{HS}_1(I))$.

Recall that $\hat{\imath}$ denotes the $i$th canonical basis vector of $\mathbb{R}^n$, and suppose that $\mathbf{1} \in \mathbb{R}^n$ denotes the vector whose all entries are equal to one. We call a monomial $u \in S$ *quasi-squarefree* if it is $\mathbf{c}$-bounded by $\mathbf{c} = \mathbf{1} + \hat{\imath}$ for some $i$. Let $I \subseteq S$ be a monomial ideal. The ideal $I$ is called quasi-squarefree monomial ideal if it is generated by quasi-squarefree monomials. For a monomial ideal $J \subseteq S$, we define its quasi-squarefree part, denoted by $J^*$, to be the ideal generated by quasi-squarefree elements of $J$.

**Lemma 2.5.** *Let $I$ be a quasi-squarefree polymatroidal ideal. Then $\mathrm{A}(I)^*$ is also a polymatroidal ideal.*

*Proof.* Suppose that $w, w' \in \mathrm{A}(I)$ are quasi-squarefree monomials in the minimal set of monomial generators of $\mathrm{A}(I)$, and suppose that $\deg_i w > \deg_i w'$ for some $i$. We show that there exists $p \in [n]$ such that $\deg_p w' > \deg_p w$ and $(w/x_i)x_p$ is a quasi-squarefree element of $\mathrm{A}(I)$. By Theorem 2.1, we know that $\mathrm{A}(I)$ is polymatroidal, and consequently, the exchange property of polymatroidal ideals guarantees the existence of $j$ such that $\deg_j w' > \deg_j w$ and $(w/x_i)x_j \in \mathrm{A}(I)$. If $(w/x_i)x_j$ is a quasi-squarefree monomial, we are done. Otherwise, assume that $(w/x_i)x_j$ is not quasi-squarefree. While $w$ is quasi-squarefree, this implies that degree of $(w/x_i)x_j$ in two distinct variables is two. Hence one has

- $\deg_j w = 1$;
- $\deg_i w = 1$;
- There exists $t$ such that $\deg_t w = 2$.

If $d(w, w') = 1$, then $(w/x_i)x_j = w'$ which contradicts our assumption that $(w/x_i)x_j$ is not quasi-squarefree. So $d(w, w') \geq 2$. On the other hand, $w$ and $w'$ are of the same degree. Hence there exists $p \neq j$ such that $\deg_p w' > \deg_p w$. Since $w'$ is quasi-squarefree, and $1 = \deg_j w < \deg_j w'$, we conclude that $\deg_p w' = 1$. Thus $\deg_p w = 0$. Let

$$\tilde{w} = (w/x_i)x_j = x_j^2 x_t^2 x_{\ell_1} \dots x_{\ell_k}$$

such that $j, t, \ell_1, \dots, \ell_k$ are pairwise ditinct. Since $\tilde{w}$ belongs to $\mathrm{A}(I)$, there exist adjacent monomials $u, u' \in \mathrm{G}(I)$ such that $\tilde{w} = \mathrm{lcm}(u, u')$. Regarding the fact that $u$ and $u'$ are quasi-squarefree, we set

$$u = x_j x_t^2 x_{\ell_1} \dots x_{\ell_k}$$

and

$$u' = x_j^2 x_t x_{\ell_1} \dots x_{\ell_k}$$

Considering the symmetric exchange property of polymatroidal ideals for $u \in \mathrm{G}(I)$ with $\deg_p u \leq \deg_p w = 0$ and an element of $\mathrm{G}(I)$ divisible by $x_p$, one deduces that $(u/x_r)x_p \in \mathrm{G}(I)$ for some $r \in \{j, t, \ell_1, \dots, \ell_k\}$. On the other hand, $(u/x_r)x_p$ is adjacent to $u$. Moreover,

$$\mathrm{lcm}((u/x_r)x_p, u) = ux_p = (w/x_i)x_p.$$

Hence the quasi-squarefree monomial $(w/x_i)x_p$ belongs to $\mathrm{A}(I)$, as desired. □

Notice that while the adjacency ideal of a polymatroidal ideal is polymatroidal as well by Theorem 2.1, in Lemma 2.5 we can not replace $\mathrm{A}(I)$ with the ideal $I$ itself; as the following example shows.

**Example 2.6.** Consider the polymatroidal ideal

$$I = (xyz^2, xyzw, y^2z^2, y^2zw) \subseteq k[x, y, z, w].$$

One can see that its quasi-squarefree part

$$I^* = (xyz^2, xyzw, y^2zw)$$

is not polymatroidal anymore. However, $\mathrm{A}(I)$ and the quasi-squarefree part of $\mathrm{A}(I)$

$$\mathrm{A}(I)^* = (xy^2zw, xyz^2w)$$

are polymatroidal.

The next lemma states that by taking iterated adjacency ideals and then quasi-squarefree part, one after another, one can obtain homological shift ideals of a polymatroidal ideal generated in degree 2.

**Lemma 2.7.** *Let $I$ be a polymatroidal ideal generated in degree two. Then $\mathrm{HS}_i(I)$ is obtained by taking successively $i$ times iterated adjacency ideals and then quasi-squarefree part starting from $I$.*

*Proof.* Assume that the elements of $I$ are ordered decreasingly in the lexicographical order with respect to $x_1 > \cdots > x_n$. By [12, Theorem 1.3], the ideal $I$ has linear quotients with respect to this ordering of generators. We first show that $\mathrm{HS}_{i+1}(I)$ is a subset of $\mathrm{A}(\mathrm{HS}_i(I))^*$, that is, the quasi-squarefree part of the adjacency ideal of $\mathrm{HS}_i(I)$. Let $ux_{\ell_1}\dots x_{\ell_{i+1}}$ be an element of $G(\mathrm{HS}_{i+1}(I))$ with $u \in \mathrm{G}(I)$ and

$$x_{\ell_1}, \dots, x_{\ell_{i+1}} \in \mathrm{set}(u).$$

Notice that if $u = x_j^2$ for some $j \in [n]$, then $x_j \notin \mathrm{set}(u)$. Hence $ux_{\ell_1}\dots x_{\ell_{i+1}}$ is a quasi-squarefree monomial. On the other hand, this element is the least common multiple of the adjacent monomials $ux_{\ell_1}\dots x_{\ell_i}$ and $ux_{\ell_1}\dots x_{\ell_{i-1}}x_{\ell_{i+1}}$ in $\mathrm{G}(\mathrm{HS}_i(I))$. So

$$\mathrm{HS}_{i+1}(I) \subseteq \mathrm{A}(\mathrm{HS}_i(I))^*.$$

Next, we show that

$$\mathrm{A}(\mathrm{HS}_i(I))^* \subseteq \mathrm{HS}_{i+1}(I). \tag{9}$$

For this purpose, let $w$ and $w'$ be two adjacent elements of $\mathrm{HS}_i(I)$ for which $\mathrm{lcm}(w, w')$ is a quasi-squarefree monomial. We are going to show that $\mathrm{lcm}(w, w') \in \mathrm{HS}_{i+1}(I)$. If $w$ and $w'$ are both squarefree monomials, then $w, w' \in \mathrm{HS}_i(I)^{\leq \mathbf{1}}$ where $\mathbf{1} = (1, \dots, 1) \in \mathbb{Z}^n$; see Section 1 for the definition of $\mathbf{1}$-bounded part $\mathrm{HS}_i(I)^{\leq \mathbf{1}}$ of $\mathrm{HS}_i(I)$. But $\mathrm{HS}_i(I)^{\leq \mathbf{1}} = \mathrm{HS}_i(I^{\leq \mathbf{1}})$ by [9, Corollary 1.10]. On the other hand, for matroidal ideal $I^{\leq \mathbf{1}}$ by [1, Corollary 3.3] one has

$$\mathrm{HS}_{i+1}(I^{\leq \mathbf{1}}) = \mathrm{A}(\mathrm{HS}_i(I^{\leq \mathbf{1}})).$$

Hence

$$\mathrm{lcm}(w, w') \in \mathrm{A}(\mathrm{HS}_i(I^{\leq \mathbf{1}}) = \mathrm{HS}_{i+1}(I^{\leq \mathbf{1}}) \subseteq \mathrm{HS}_{i+1}(I).$$

So we are done in the case that $w$ and $w'$ are both squarefree.

Next assume that $w$ is not squarefree. Suppose that $w = ux_{\ell_1}\dots x_{\ell_i}$ with $u \in \mathrm{G}(I)$, and $x_{\ell_1}, \dots, x_{\ell_i}$ are elements of $\mathrm{set}(u)$. We are going to show that $\mathrm{lcm}(w, w')$ can be written as a product of an element of $\mathrm{G}(I)$ and $i + 1$ distinct elements of its set which implies the desired inclusion (9). Since $w$ and $w'$ are of the same degree and adjacent, the monomial $w' : w$ is of degree one, say

$$w' : w = \frac{w'}{\gcd(w, w')} = x_f. \tag{10}$$

Hence

$$\operatorname{lcm}(w, w') = wx_f.$$

Regarding our assumption that $w$ is not squarefree, there exists $j \in [n]$ such that $\deg_j w = 2$. There are two cases to consider:

Case 1. Let $u = x_j^2$. Besides, there exist monomials in $\mathrm{G}(I)$ divided by each variable $x_{\ell_1}, \dots, x_{\ell_i}$ and $x_f$; a pairwise distinct sequence of variables as discussed below. Applying the exchange property of polymatroidal ideals for $u = x_j^2$ and the monomials divisible by these mentioned variables, one concludes that $x_j x_{\ell_1}, \dots, x_j x_{\ell_i}$ and $x_j x_f$ are elements of $\mathrm{G}(I)$. Now consider the following elements of $\mathrm{G}(I)$:

$$x_j^2, x_j x_{\ell_1}, \dots, x_j x_{\ell_i}, x_j x_f. \tag{11}$$

These elements are pairwise distinct because

- $x_{\ell_1}, \dots, x_{\ell_i} \in \operatorname{set}(u)$ are pairwise distinct by Remark 1.1;
- the ideal $I$ is generated in degree two. This implies that $x_j \notin \operatorname{set}(x_j^2)$. Consequently, $x_j \notin \{x_{\ell_1}, \dots, x_{\ell_i}\}$;
- $x_j \neq x_f$. Otherwise, regarding (10), $\deg_f w' = 3$ which contradicts the fact that $w'$ is quasi-squarefree;
- $x_f \notin \{x_{\ell_1}, \dots, x_{\ell_i}\}$. Otherwise, $\deg_f w = 1$ and consequently, $\deg_f w' = 2$ by (10). So $x_f^2$ and $x_j^2$ both divide $\operatorname{lcm}(w, w')$, a contradiction to the assumption that $\operatorname{lcm}(w, w')$ is quasi-squarefree.

Let $v$ be the maximum element in (11) with respect to the lexicographical order induced by $x_1 > \cdots > x_n$. By considering $\tilde{v} : v$ for each

$$\tilde{v} \in \{x_j^2, x_j x_{\ell_1}, \dots, x_j x_{\ell_i}, x_j x_f\} \setminus \{v\},$$

one obtains $i+1$ distinct varibles $x_{k_1}, \dots, x_{k_{i+1}}$ in $\operatorname{set}(v)$. One can see that

$$v x_{k_1} \dots x_{k_{i+1}} = x_j^2 x_{\ell_1} \dots x_{\ell_i} x_f = \operatorname{lcm}(w, w').$$

Thus $\operatorname{lcm}(w, w') \in \mathrm{HS}_{i+1}(I)$.

Case 2. Let $u = x_j x_p$ for some $p \neq j$, and $x_j \in \{x_{\ell_1}, \dots, x_{\ell_i}\}$. We first discuss why in this case $x_{\ell_1}, \dots, x_{\ell_i}, x_p, x_f$ are pairwise distinct. For this purpose, we notice that

- $x_{\ell_1}, \dots, x_{\ell_i} \in \operatorname{set}(u)$ are pairwise distinct as we have seen in Case 1.
- $x_p \notin \{x_{\ell_1}, \dots, x_{\ell_i}\}$. Otherwise, $x_p^2$ and $x_j^2$ both divide $w$, a contradiction to the fact that $w$ is quasi-squarefree.
- $x_f \notin \{x_{\ell_1}, \dots, x_{\ell_i}\}$. Otherwise, $\deg_f w = 1$, and consequently $\deg_f w' = 2$ by (10). This implies that $\operatorname{lcm}(w, w')$ is divisible by $x_f^2$ and $x_j^2$ while $j \neq f$ regarding the degree of $w$ and $w'$ in these variables. As a result, $\operatorname{lcm}(w, w')$ is not quasi-squarefree, a contradiction.
- $x_f \neq x_p$. Otherwise, $\deg_f w' = 1 + \deg_f w = 2$. On the other hand, $f \neq j$ as discussed above. As a result distinct monomials $x_j^2$ and $x_f^2$ divide $\operatorname{lcm}(w, w')$ which contradicts its property of being quasi-squarefree.

Based on this discussion, $x_{\ell_1}, \dots, x_{\ell_i}, x_p, x_f$ are pairwise distinct. Since $x_j \in \{x_{\ell_1}, \dots, x_{\ell_i}\} \subseteq \text{set}(x_j x_p)$, we deduce that $x_j^2 \in I$. On the other hand, for each variable $x_{\ell_1}, \dots, x_{\ell_i}$ and $x_f$, there exists a monomial which is divisible by that variable. Applying the exchange property of polymatroidal ideals for $x_j^2$ and those monomials, we conclude that $x_j x_{\ell_1}, \dots, x_j x_{\ell_i}, x_j x_f \in \mathrm{G}(I)$. Recall that $u = x_j x_p$ also belongs to $\mathrm{G}(I)$. Furthermore,

$$x_j x_{\ell_1}, \dots, x_j x_{\ell_i}, x_j x_f, x_j x_p$$

are pairwise distinct as discussed above. Suppose that $v$ is the maximum element of $\{x_j x_{\ell_1}, \dots, x_j x_{\ell_i}, x_j x_f, x_j x_p\}$ with respect to the lexicographical order induced by $x_1 > \cdots > x_n$. Now by considering

$$\tilde{v} : v \quad \text{for each} \quad \tilde{v} \in \{x_j x_{\ell_1}, \dots, x_j x_{\ell_i}, x_j x_f, x_j x_p\} \setminus \{v\},$$

we obtain $i+1$ distinct variables $x_{k_1}, \dots, x_{k_{i+1}}$ in $\text{set}(v)$. Moreover,

$$v x_{k_1} \dots x_{k_{i+1}} = (x_j x_p) x_{\ell_1} \dots x_{\ell_i} x_f = \text{lcm}(w, w'),$$

as desired.

**Theorem 2.8.** *Let $I$ be a polymatroidal ideal generated in degree two. Then*

$$\mathrm{HS}_{i+j}(I) \subseteq \mathrm{HS}_i(\mathrm{HS}_j(I))$$

*for each $i, j$.*

*Proof.* It is enough to notice that by Lemma 2.5 and Lemma 2.7, the ideal $\mathrm{HS}_j(I)$ is polymatroidal for each $j$. Thus the assertion follows from Theorem 1.2. □

As an immediate consequence of Lemma 2.5 and Lemma 2.7 we have

**Corollary 2.9.** ([7, Theorem 4.5]) *Let $I$ be a polymatroidal ideal generated in degree two. Then the ideal $\mathrm{HS}_i(I)$ is also a polymatroidal ideal for all $i$.*

**Proposition 2.10.** ([5, Proposition 2.4]) *Let $I$ be a matroidal ideal. Then*

$$\mathrm{HS}_{i+j}(I) = \mathrm{HS}_i(\mathrm{HS}_j(I))$$

*for each $i, j$.*

*Proof.* By [1, Corollary 3.3], the ideal $\mathrm{HS}_{i+j}(I)$ can be obtained by taking $i+j$ times iterated adjacency ideals starting from $I$ for each $i, j$. On the other hand, $\mathrm{HS}_j(I)$ is also a matroidal ideal by [1, Theorem 3.2]. So one obtains $\mathrm{HS}_i(\mathrm{HS}_j(I))$ by taking first $j$ times, and next $i$ more times adjacency ideals starting from $I$. □

# References

[1] S. Bayati, *Multigraded shifts of matroidal ideals*, Archiv der Mathematik 111 (3), 239–246 (2018)

[2] S. Bayati, I. Jahani, N. Taghipour, *Linear quotients and multigraded shifts of Borel ideals*, Bulletin of the Australian Mathematical Society 100, no. 1, 48–57 (2019)

[3] A. Conca, J. Herzog, *Castelnuovo-–Mumford regularity of products of ideals*, Collect. Math. 54, 137–152 (2003)

[4] S. Eliahou, M. Kervaire, *Minimal resolution of some monomial ideals*, Journal of Algebra 129, no. 1, 1–25 (1990)

[5] A. Ficarra, *Homological shifts of polymatroidal ideals*, arXiv:2205.04163v2 [math.AC] (2022)

[6] V. Gasharova, T. Hibi, I. Peeva, *Resolutions of a-stable ideals*, Journal of Algebra 254 (2), 375–394 (2002)

[7] A. Ficarra, J. Herzog, *Dirac's theorem and multigraded syzygies*, arXiv: 2212.00395v1 [math.AC] (2022)

[8] J. Herzog, T. Hibi, X. Zheng, *Dirac's theorem on chordal graphs and Alexander duality*, European Journal of Combinatorics 25 (7), 949–960 (2004)

[9] J. Herzog, S. Moradi, M. Rahimbeigi, G. Zhu, *Homological shift ideals*, Collectanea Mathematica 72, 157–174 (2021)

[10] J. Herzog, S. Moradi, M. Rahimbeigi, G. Zhu, *Some homological properties of Borel type ideals*, arXiv:2112.11726v1 [math.AC] (2021)

[11] J. Herzog, Y. Takayama, *Resolutions by mapping cones*, The Roos Festschrift volume 4, Homology, Homotopy and Applications 2 (2), 277–294 (2002)

[12] F. Mohammadi, S. Moradi, *Weakly polymatroidal ideals with applications to vertex cover ideals*, Osaka J. Math. 47, no. 3, 627–-636 (2010)

Department of Mathematics and Computer Science, Amirkabir University of Technology (Tehran Polytechnic), Iran

*Email address*: shamilabayati@gmail.com